\documentclass[12pt]{amsart}

\usepackage{tgpagella}
\usepackage{euler}
\usepackage[T1]{fontenc}
\usepackage{amsmath, amssymb, mathabx}
\usepackage[hidelinks]{hyperref}
\usepackage[english]{babel}
\usepackage{mathrsfs}
\usepackage{eucal}
\usepackage[all]{xy}
\usepackage{tikz}

\newtheorem{thm}{Theorem}[section]
\newtheorem*{thm*}{Theorem}
\newtheorem{lem}[thm]{Lemma}

\newtheorem{prop}[thm]{Proposition}
\newtheorem*{prop*}{Proposition}

\newtheorem{cor}[thm]{Corollary}
\newtheorem*{cor*}{Corollary}

\theoremstyle{definition}
\newtheorem{defn}[thm]{Definition}
\newtheorem*{defn*}{Definition}

\newtheorem{question}[thm]{Question}

\newtheorem*{question*}{Question}
\newtheorem*{Pquestion*}{Popa's question}

\newtheorem*{conv*}{Convention}

\newcommand{\dminus}{ 
\buildrel\textstyle\ .\over{\hbox{ 
\vrule height3pt depth0pt width0pt}{\smash-} 
}}

\def\bb{\mathbb}

\def\bb{\mathbb}

\def\cal{\mathcal}

\def\u{\mathsf 1}

\newcommand{\cstar}{$\mathrm{C}^*$}

\makeatletter

\def\dotminussym#1#2{%
  \setbox0=\hbox{$\m@th#1-$}%
  \kern.5\wd0%
  \hbox to 0pt{\hss\hbox{$\m@th#1-$}\hss}%
  \raise.6\ht0\hbox to 0pt{\hss$\m@th#1.$\hss}%
  \kern.5\wd0}

\newcommand\ab[1]{\left| #1 \right|}

\def \Th{\operatorname{Th}}
\def \R{\mathcal R}
\def \u{\mathcal U}

\def \val{\operatorname{val}}
\def \pval{\operatorname{s-val}}

\def\O{\mathcal{O}}

\def \zdcomp{\operatorname{ZDComp}}
\def \comp{\operatorname{Comp}}
\def \bool{\operatorname{Bool}}
\def \ab{\operatorname{AbC*}}
\def \rrz{\operatorname{RRZ}}
\def \nuc{\operatorname{nuc}}

\allowdisplaybreaks
\makeatletter
\def\l@subsection{\@tocline{2}{0pt}{2.5pc}{5pc}{}}
\def\l@subsubsection{\@tocline{2}{0pt}{5pc}{7.5pc}{}}
\makeatother

\begin{document}

%%%%%%%%%%%%%%%%%%%%%%%%%%%%%%%%%%%%%%%%%%%%%%

\title{Model theory and Ultrapower Embedding Problems in Operator Algebras}

\author{Isaac Goldbring}
\address{Department of Mathematics\\University of California, Irvine, 340 Rowland Hall (Bldg.\# 400),
Irvine, CA 92697-3875}
\email{isaac@math.uci.edu}
\urladdr{http://www.math.uci.edu/~isaac}
\thanks{Goldbring was partially supported by NSF grant DMS-2054477.}

\begin{abstract}
We survey the model theoretic approach to a variety of ultrapower embedding problems in operator algebras. 
\end{abstract}

\maketitle

% \begin{abstract}

% \end{abstract}

\section{Introduction}

Suppose that $T$ is an $\forall\exists$-axiomatizable theory in some countable, classical language $L$.  It is often too much to ask that there exists a countable model of $T$ into which all other countable models embed.  (This does happen, for example, when the theory $T$ admits an $\aleph_0$-categorical model companion.)  Nevertheless, under the mild assumption that the theory $T$ has the \emph{joint embedding property} (JEP), meaning that any two models of $T$ mutually embed into a third model of $T$, we can infer the existence of countable models of $T$ whose ultrapower with resepct to any nonprincipal ultrafilter $\u$ on $\bb N$ embeds all countable models of $T$; we refer to such models of $T$ as \emph{locally universal} models of $T$.  Indeed, if $M$ is any \emph{existentially closed} model of $T$ and $N$ is any countable model of $T$, then by jointly embedding $M$ and $N$ into a (without loss of generality) countable model $P$ of $T$, we see that $P$ embeds in $M^\u$ (this follows from the fact that $M$ is e.c. and $M^\u$ is a countably saturated model of its theory) whence $N$ also embeds into $M^\u$.  (By using so-called \emph{good ultrafilters}, one can obtain ultrapowers of $M$ which embed larger models of $T$.)  This discussion holds verbatim for continuous theories:  if $T$ is an $\forall\exists$-axiomatizable theory with JEP in a separable language, then separable locally universal models of $T$ exist.  It is clear that any model of $T$ that contains a locally universal model of $T$ is itself locally universal, whence countable (separable) locally universal models of $T$ are ubiquitous.

Seemingly unaware of this abstract model-theoretic discussion, operator algebraists have posed a number of problems which ask whether or not concrete operator algebras are locally universal for the corresponding classes to which they belong (which are not always elementary classes).  The most famous of these problems is the \emph{Connes embedding problem}, which appeared in Alain Connes' fundamental work \cite{connes} from 1976, in which he showed that any separable injective II$_1$ factor is necessarily hyperfinite (and for which he received the Fields medal in 1982).  Part of the proof of his main theorem involved showing that a particular separable II$_1$ factor embedded into an ultrapower of the \emph{hyperfinite II$_1$ factor} $\R$.  He casually remarked that such an embedding ``ought to'' exist for any separable II$_1$ factor, that is, $\R$ should be a locally universal object for the elementary class of II$_1$ factors.  (Incidentally, since $\R$ embeds into any II$_1$ factor, the CEP is equivalent to the assertion that all II$_1$ factors have the same universal theory.)  Connes' reason why this ``ought to'' be the case is not entirely compelling:  he points out that such an embedding exists for $L(\bb F_2)$, the group von Neumann algebra associated to the free group on two generators, and it is for this reason that such an embedding should exist for all separable II$_1$ factors.  While some operator algebraists refer to Connes' offhand remark as a ``Conjecture,'' most prefer to call it a ``Problem.''

Over the years, the CEP has gained significant interest due its connections with a wide variety of areas of mathematics, including  \cstar-algebra theory (in connection with \emph{Kirchberg's QWEP problem}), quantum information theory (in connection with \emph{Tsirelson's problem}), free probability, group theory, and noncommutative real algebraic geometry, to name a few.  

Remarkably, in early 2020, a negative resolution was obtained to the CEP via its equivalence with Tsirelson's problem, which was itself refuted using a remarkable theorem in quantum complexity theory known as $\operatorname{MIP}^*=\operatorname{RE}$.  Independent of its connection with the CEP, this latter result is widely considered to be a landmark scientific achievement; the reader interested in understanding the entire story behind these connections can consult the author's survey \cite{bams}.

That being said, someone wishing to understand the proof of the negative solution to the CEP using $\operatorname{MIP}^*=\operatorname{RE}$ must tread the deep waters connecting these two seemingly distant results.  In joint work with Bradd Hart \cite{undec}, we showed how, using basic ideas from continuous model theory (most notably the Completeness Theorem for first-order continuous logic and the theory of definable sets in continuous logic), one can obtain a more direct proof of the failure of CEP from $\operatorname{MIP}^*=\operatorname{RE}$.  Moreover, the model-theoretic approach offers more insight into this refutation and allows one to prove extra results, such as ``many counterexamples'' to the CEP, that is, many different universal theories of II$_1$ factors, as well as a  G\"odelian-style refutation stating that any effectively axiomatizable class of II$_1$ factors will contain a counterexample to the CEP.  We present this model theoretic approach to the negative solution to the CEP in Section \ref{CEPsection}.

In Sections \ref{kepsection} and \ref{MFsection}, we consider two \cstar-algebraic analogs of the CEP, the so-called \emph{Kirchberg embedding problem} and the \emph{MF problem}, which, in some sense, can be thought of as the ``infinite'' and ``finite'' \cstar-versions of the CEP.  The former problem asks whether or not the \emph{Cuntz algebra} $\cal O_2$ is locally universal for the class of all \cstar-algebras.  This problem has eluded model-theoretic techniques thus far and we discuss what we know about this still-open problem; the material presented in this section represents joint work of the author and Thomas Sinclair \cite{GS}.  The latter problem is the so-called \emph{MF problem}, which asks whether or not the \emph{universal UHF algebra} $\cal Q$ is a locally universal object for the elementary class of stably finite \cstar-algebras.  While a negative solution to the MF problem can readily be deduced from the failure of the CEP, the model-theoretic approach allows us to establish a similar G\"odelian-style refutation of the MF problem.  

A variant of the MF problem, known as the \emph{quasidiagonality (QD) problem}, asks if $\cal Q$ is a locally universal object for the (non-elementary!) class of stably finite \emph{nuclear} \cstar-algebras.  While the QD problem is still open, a major breakthrough was obtained by Tikuisis, Winter, and White \cite{TWW}, who resolved this problem in the positive (for simple such algebras) assuming a technical assumption known as the \emph{universal coefficient theorem} (UCT).  The model-theoretic content of the UCT is not widely understood at the time of the writing of this paper.  Nevertheless, we discuss some model-theoretic ideas around the QD problem representing more joint work of the author and  Sinclair \cite{international}.  

In the final subsection of Section \ref{kepsection}, we return to the ideas of the author and Hart from \cite{undec} and show how the G\"odelian-style refutation of the MF problem extends to a larger class of algebras which, in particular, allow us to refute a \emph{stably projectionless} version of the MF problem which asks if every stably projectionless algebra embeds into an ultrapower of a very important algebra in the classification program for nuclear \cstar-algebras known as the \emph{Jiang Su algebra} $\cal Z$.  Currently, this latter result has no purely operator-algebraic proof.

In the final section, we consider the simpler case of (unital) abelian \cstar-algebras.  There, an $\aleph_0$-categorical model completion exists, namely the theory of $C(2^{\bb N})$. This result is essentially (after some category-theoretic considerations) a restatement of the existence of an $\aleph_0$-categorical model completion for the classical theory of Boolean algebras, namely the theory of atomless Boolean algebras.  However, an interesting phenomenon arises when restricting to the class of \emph{projectionless} abelian \cstar-algebas, whose models are of the form $C(X)$ for $X$ a \emph{continuum} (that is, a connected compact Hausdorff space).  In this case, a theorem of K.P. Hart \cite{KP} states that all such algebras (except for the trivial case of $C(\text{a point})\cong \bb C$) have the same universal theory, whence all nontrivial objects of this class are locally universal!  We provide a fairly detailed proof of this result below.  However, this theorem does not represent the end of the story for this class of algebras, for the question of a model companion for this class is still open.  Concerning existentially closed projectionless abelian \cstar-algebras, we discuss the positive solution, due to Christopher Eagle, Alessandro Vignati, and the author \cite{EGV} of a question of Bankston, who asked if a fairly important (and generic) continuum, the so-called \emph{pseudoarc} $\bb P$, is co-existentially closed (which just means that $C(\bb P)$ is an existentially closed projectionless abelian \cstar-algebra).

Many of the results to be discussed below involve a particular approach to the Henkin construction known as building models by games (as first popularized by Hodges in his book \cite{hodges}).  We discuss the essential properties of this construction in the next section.

\section{Building models by games}\label{gamessection}

In this section, we recall the basic facts from the theory of model-theoretic forcing needed throughout this paper.  The relevant version of model-theoretic forcing for us is the game-theoretic approach, originally presented in Hodges' classic book \cite{hodges} and adapted to the continuous setting by the author in \cite{enforceable}.  That being said, for some of what is to follows, we need to consider a slightly more general setting and so we take the opportunity to extend the context here.

Throughout this section, we fix a countable (continuous) language $L$.  By a \textbf{$\forall\bigvee\exists$-sentence} we mean an $L_{\omega_1,\omega}$-sentence of the form $$\sup_{x_1}\cdots\sup_{x_{k_n}}\bigvee_{m\in \bb N}\varphi_m(x_1,\ldots,x_{k_n}),$$ where each $\varphi_m$ is an existential $L$-formula and the symbol $\bigvee$ denotes a countable infimum.  (We note that there are several different approaches to $L_{\omega_1,\omega}$ in the literature; in the above sentences, there is no requirement on a common modulus of uniform continuity for the formulae appearing in the countable infimum).  By a \textbf{$\forall\bigvee\exists$-theory} we mean a collection of $\forall\bigvee\exists$-sentences.  We say that a class $\cal K$ of $L$-structures is \textbf{$\forall\bigvee\exists$-axiomatizable} if there is a $\forall\bigvee\exists$-theory $T$ such that, for all $L$-structures $A$, we have $A\in \cal K$ if and only if $\sigma^A=0$ for all $\sigma\in T$.  The need to consider such infinitary theories arises as many important classes of \cstar-algebras are not first-order axiomatizable but are $\forall\bigvee\exists$-axiomatizable, such as simple \cstar-algebras and nuclear \cstar-algebras (see, for example, \cite{munster}, where the class is called \emph{definable by a uniform family of formulae}).  Motivated by this, if $T$ is a $\forall\bigvee\exists$-axiomatizable theory and $P$ is a property that may or may not hold of models of $T$, we say that $P$ is $\forall\bigvee\exists$-axiomatizable (relative to $T$) if the collection of models of $T$ having property $P$ is itself $\forall\bigvee\exists$-axiomatizable.

In the remainder of this subsection, we fix a $\forall\bigvee\exists$-axiomatizable $L$-theory $T$.  We note that this assumption implies that the class of models of $T$ is an inductive class (that is, is closed under direct limits), whence every model of $T$ is contained in an e.c. model of $T$ of the same density character.

We now fix a countably infinite set $C$ of constant symbols enumerated $(c_n)_{n<\omega}$.  A \textbf{condition} (relative to $T$) is a finite set $p$ of expressions of the form $\varphi<r$, where $\varphi$ is a quantifier-free $L(C)$-sentence, such that $T\cup p$ is satisfiable.

We now consider a two-player game whose players we denote by $\forall$ and $\exists$.  The players take turns playing conditions with the requirement that each player's turn extends the previous player's turn.  When they have finished the game, they have constructed a countably infinite set $\bar p=\bigcup_{n<\omega} p_n$ which is the union of all of the conditions played during the game.  We call the play of the game \textbf{definitive} if:  for every atomic $L(C)$-sentence $\theta$, there is a unique $r\in \bb R$ such that $T\cup \bar p\models |\theta-r|=0$.  In this case, $\bar p$ uniquely describes a separable $L(C)$-structure $A^+(\bar p)$ generated by $C$, called the \textbf{compiled structure}, whose $L$-reduct is denoted $A(\bar p)$.  Since $\exists$ can always ensure that the play of the game is definitive, we assume from here on out that all plays of all games are definitive, whence determine compiled structures.

If $P$ is a property of $L(C)$-structures, we say that $P$ is \textbf{enforceable} if $\exists$ has a strategy that ensures that the compiled structure has property $P$.  It is routine to check that a countable conjunction of enforceable properties is once again enforceable.

For us, one of the most important facts about enforceability is that the property ``$A(\bar p)$ is an e.c. model of $T$'' is enforceable.  While proven in \cite{enforceable} for the case of an $\forall\exists$-theory $T$, the proof readily extends to this more general case.

We will also need the following fact:

\begin{prop}\label{force}
If $P$ is a $\forall\bigvee\exists$-axiomatizable property (relative to $T$) and there is a locally universal model of $T$ with property $P$, then $P$ is an enforceable property.
\end{prop}

Once again, this was proven in the first-order context in \cite{enforceable}, but the proof readily extends.  This construction of e.c. models of $T_\forall$ with $\forall\bigvee\exists$-axiomatizable properties is sometimes called the ``Eastern form'' of the omitting types theorem (see \cite{keislerhandbook}).

Finally, we say that an $L$-structure $A$ is an \textbf{enforceable structure} (or an \textbf{enforceable model of $T$}) if the property ``$A(\bar p)\cong A$'' is an enforceable property.  It is clear that there can be at most one enforceable structure (up to isomorphism), which is then necessarily an e.c. model of $T$.

\section{The Connes Embedding Problem}\label{CEPsection}

Recall from the introduction that the Connes Embedding Problem (CEP) asks whether or not the hyperfinite II$_1$ factor $\R$ is a locally universal model of the theory of tracial von Neumann algebras.  In this section, we discuss the recent resolution of this problem in the negative.  A much more detailed version of this section appears in the author's survey article \cite{bams}.

% As mentioned in the introduction, the Connes Embedding Problem (that is, the statement that $\R$ is a locally universal II$_1$ factor) was refuted by a statement in quantum complexity theory known as $\operatorname{MIP}^*=\operatorname{RE}$.  We begin by explaining the version of this theorem that is easiest to explain and briefly describe how this result led to the negative solution to the CEP.  We then discuss the model-theoretic explanation for this implication.

\subsection{Nonlocal games and their entangled strategies}  

The negative solution to the CEP follows from a result in quantum complexity theory known as $\operatorname{MIP}^*=\operatorname{RE}$, which we now proceed to explain.  The central object at the heart of this result is the notion of a nonlocal game.

\begin{defn}
A \textbf{nonlocal game with $k$ questions and $n$ answers} is a pair $\frak G:=(\pi,D)$, where $\pi$ is a probability distribution on $[k]\times [k]$ and $$D:[k]\times [k]\times [n]\times [n]\to \{0,1\}$$ is a function, known as the \textbf{decision predicate}.  We refer to the pair $(k,n)$ as the \textbf{dimensions of the game} (although this terminology is not standard).
\end{defn}

Here, $[k]:=\{1,\ldots,k\}$ and similarly for $[n]$.  One envisions the game $\frak G$ above being played by two players, usually named Alice and Bob.  The referee for the game randomly sends Alice and Bob a pair of questions $(x,y)\in [k]\times [k]$ according to the distribution $\pi$.  Somehow, Alice and Bob respond with answers $a$ and $b$ respectively from $[n]$.  The referee then uses the function $D$ to decide if they won this particular play of the game, that is, they win if and only if $D(x,y,a,b)=1$.

In the result $\operatorname{MIP}^*=\operatorname{RE}$, Alice and Bob use so-called entangled strategies to decide how to play:

\begin{defn}

\

\begin{enumerate}
    \item If $\cal H$ is a Hilbert space, a \textbf{positive operator-valued measure} (or \textbf{POVM} for short) \textbf{of length $n$} is a sequence $A_1,\ldots,A_n$ of positive operators on $\cal H$ such that $\sum_{i=1}^n A_i=I_{\cal H}$.
    \item An \textbf{entangled strategy} for $(k,n)$-dimensional games consists of the following information:
    \begin{enumerate}
        \item A pair of finite-dimensional Hilbert spaces $\cal H_A$ and $\cal H_B$;
        \item For each $x\in [k]$, a POVM $A^x_1,\ldots,A^x_n$ on $\cal H_A$;
        \item For each $x\in [k]$, a POVM $B^x_1,\ldots,B^x_n$ on $\cal H_B$;
        \item A unit vector $\psi\in \cal H_A\otimes \cal H_B$.
    \end{enumerate}
\end{enumerate}
\end{defn}

The above definition represents the situation that Alice and Bob will perform quantum experiments to determine how to answer the questions they are sent by the referee.  Alice has a quantum system with associated Hilbert space $\cal H_A$ while Bob has another quantum system with associated Hilbert space $\cal H_B$.  They share some state $\psi$ of their composite system $\cal H_A\otimes \cal H_B$, which is usually an entangled state (whence the name entangled strategy), that is, is not necessarily a simple tensor $\psi_A\otimes \psi_B$ with $\psi_A\in \cal H_A$ and $\psi_B\in \cal H_B$.  Upon receiving question $x$, Alice performs the measurement corresponding to the POVM $A^x$ on her part of the state $\psi$ to determine how she will respond; Bob acts in a similar fashion.

Given an entangled strategy for $(k,n)$-dimensional games as above, we calculate, for $x,y\in [k]$ and $a,b\in [n]$, the value 
$$p(a,b|x,y):=\langle (A^x_a\otimes B^y_b)\psi,\psi\rangle,$$ which we interpret as the probability that Alice and Bob respond with answers $a$ and $b$ if asked questions $x$ and $y$ respectively.  If we let $p$ denote the strategy thus obtained, then for any nonlocal game $\frak G$ with $k$ questions and $n$ answers, the \textbf{entangled value of $\frak G$ corresponding to $p$} is defined to be 
$$\val(\frak G,p):=\sum_{(x,y)\in [k]\times [k]}\pi(x,y)\sum_{(a,b)\in [n]\times [n]}p(a,b|x,y)D(x,y,a,b),$$ which measures the probability that Alice and Bob win the game $\frak G$ if they play according to the strategy $p$.  We let $C_q(k,n)\subseteq [0,1]^{k^2n^2}$ denote the set of entangled strategies for $(k,n)$-dimensional games.  The optimal probability Alice and Bob have for winning the game by using entangled strategies is thus
$$\val^*(\frak G):=\sup_{p\in C_q(k,n)}\val(\frak G,p)$$ and is called the \textbf{entangled value of $\frak G$}.

We can now state the version of $\operatorname{MIP}^*=\operatorname{RE}$ relevant for us:

\begin{thm}[Ji et. al. \cite{MIP*}]
There is an ``efficient'' mapping $\bf M\mapsto \frak G_{\bf M}$ from Turing machines to nonlocal games such that:
\begin{itemize}
    \item If $\bf M$ halts on the empty tape, then $\val^*(\frak G_{\bf M})=1$.
    \item If $\bf M$ does not halt on the empty tape, then $\val^*(\frak G_{\bf M})\leq \frac{1}{2}$.
\end{itemize}
\end{thm}

The reader interested in the completxity-theoretic motivation (and nomenclature) for this result can consult the aforementioned survey \cite{bams}.

\subsection{From $\operatorname{MIP}^*=\operatorname{RE}$ to the failure of CEP:  the traditional route}

We now briefly describe how the above result was used to refute the CEP.  We begin by noting that one can effectively approximate $\val^*(\frak G)$ from below (provided the distribution $\pi$ can be effectively calculated).  Indeed, given some dimension $d$, one can enumerate a computable sequence of finite nets $$N^d_1\subseteq N^d_2\subseteq \cdots$$ over all states and POVMs in dimension $d$ with $|N^d_m|=m^{O(d^2)}$ such that, for any $p\in C_q(k,n)$ based on a $d$-dimensional strategy and any $m$, there is $q\in N^d_m$ with $|\val(\frak G,p)-\val(\frak G,q)|<\frac{1}{m}$.  If one sets $\val^n(\frak G,p)=\max_{d,m\leq n}\max_{p\in N^d_m}\val(\frak G,p)$,  then $\val^n(\frak G,p)$ is computable and $\val^n(\frak G,p)$ converges to $\val(\frak G)$ from below as $n$ tends to $\infty$.

By the equality $\operatorname{MIP}^*=\operatorname{RE}$, we deduce that there cannot exist an effective procedure, uniform over all games, for approximating $\val^*(\frak G)$ from above.  Indeed, if this were possible, then, combining this procedure with the procedure from the previous paragraph, one would be able to effectively approximate $\val^*(\frak G)$ uniformly in the description of $\frak G$.  Any estimate of $\val^*(\frak G_{\bf M})$ to within $\frac{1}{4}$ would allow one to decide whether or not $\val^*(\frak G_{\bf M})$ was $1$ or at most $\frac{1}{2}$, and thus one could effectively solve the halting problem!

The relevance of this conclusion is that it refutes a well-known problem in quantum information theory known as \textbf{Tsirelson's problem}.  Tsirelson's problem questioned whether or not the entangled value $\val^*(\frak G)$ of $\frak G$ and its so-called \textbf{commuting value} $\val^{co}(\frak G)$ defined using an a priori larger set of strategies known as quantum commuting strategies coincided.  Here, a \textbf{quantum commuting strategy} for $(k,n)$-dimensional games consists of:
\begin{itemize}
    \item A single (possibly infinite-dimensional) Hilbert space $\cal H$,
    \item For each $x\in [k]$, two POVMs $A^x_1,\ldots,A^x_n$ and $B^y_1,\ldots,B^y_b$ on $\cal H$ satisfying $A^x_aB^y_b=B^y_bA^x_a$ for all $x,y\in [k]\times [k]$ and $a,b\in [n]\times [n]$, and
    \item A unit vector $\psi\in \cal H$.
\end{itemize}  

The second condition in the previous definition ensures that the measurements $A^x$ and $B^y$ can be done simultaneously (ensuring that the players still cannot communicate).  Using this strategy, the probability that Alice and Bob respond with $a$ and $b$ if asked questions $x$ and $y$ is given by $p(a,b|x,y)=\langle A^x_aB^y_b\xi,\xi\rangle$.  The set of such quantum commuting strategies for $(k,n)$-dimensional games is denoted $C_{qc}(k,n)$ and the quantum commuting value of $\frak G$ is $\val^{co}(\frak G)=\sup_{p\in C_{qc}(k,n)}\val(\frak G,p)$.  It follows from the results in \cite{Fr} and \cite{FNT} that $\val^{co}(\frak G)$ can always be uniformly effectively approximated from above.  Thus, Tsirelson's problem must have a negative answer by the results of the previous paragraph.

It had already been obvserved by Fritz \cite{Fr} and independently by Junge et. al. \cite{Junge} that a positive answer to a well-known problem in \cstar-algebra theory known as \textbf{Kirchberg's QWEP problem} would yield a positive answer to Tsirelson's problem, whence Kirchberg's QWEP problem must also have a negative answer.  (Ozawa \cite{Oz} had later shown that a positive answer to Tsirelson's problem would in turn imply a positive answer to Kirchberg's problem, but that is now a moot point.)  A difficult result of Kirchberg \cite{K} states that his QWEP problem is actually equivalent to the CEP, whence we finally obtain the negative solution to the CEP from $\operatorname{MIP}^*=\operatorname{RE}$!  For details on all of this, see the aforementioned survey \cite{bams}.

\subsection{From $\operatorname{MIP}^*=\operatorname{RE}$ to the failure of CEP:  the model-theoretic route}

We now explain how continuous model theory can be used to eliminate the use of Tsirelson's problem and Kirchberg's QWEP problem in the derivation of the negative solution of CEP from $\operatorname{MIP}^*=\operatorname{RE}$.  First, we set $C_{qa}(k,n):=\overline{C_q(k,n)}$; the elements of this set are the strategies for $(k,n)$-dimensional games that can be approximated by entangled strategies.

\begin{defn}
A strategy $p\in C_{qa}(k,n)$ is said to be \textbf{synchronous} if, for all $x\in [k]$ and \emph{distinct} $a,b\in [n]$, we have $p(a,b|x,x)=0$.
\end{defn}

In other words, a strategy is synchronous if there is no chance that Alice and Bob respond differently if asked the same question.  We let $C_{qa}^s(k,n)$ denote the set of synchronous elements of $C_{qa}$ and we let $\pval(\frak G):=\sup_{p\in C_{qa}^s(k,n)}\val(\frak G,p)$ denote its \textbf{synchronous value}.  It is an artifact of the proof of $\operatorname{MIP}^*=\operatorname{RE}$ that one may replace $\val^*$ with $\pval^*$ without changing the validity of the result.  The relevance of the set $C_{qa}^s(k,n)$ for us is the following operator-algebraic characterization:

\begin{thm}[Kim, Paulsen, and Schaufhauser \cite{KPS}]\label{KPStheorem}
Given $p\in [0,1]^{k^2n^2}$, we have that $p\in C_{qa}^s(k,n)$ if and only if, for every $x\in [k]$, there are projections $e^x_1,\ldots,e^x_n\in \cal R^\u$ such that $\sum_{a=1}^n e^x_a=1$ and such that $p(a,b|x,y)=\tau(e^x_ae^y_b)$ (where $\tau$ denotes the unique trace on $\cal R^\u$).
\end{thm}

Consequently, if we set $z:=(z^x_a)$ to be variables ranging over the operator norm unit ball and set $$\varphi(z):=\sum_{(x,y)\in [k]\times [k]}\pi(x,y)\sum_{(a,b)\in [n]\times [n]}\tau(z^x_az^y_b)D(x,y,a,b),$$ a quantifier-free formula in the language of tracial von Neumann algebras, then for every $p\in C_{qa}^s(k,n)$, we have $\val(\frak G,p)=\varphi(e)^{\R^\u}$, where $e=(e^x_a)$ is the family of projections guaranteed to exist by Theorem \ref{KPStheorem}.

For any $k$ and $n$, set $Z(k,n)$ to be the set consisting of all tuples $e=(e^x_a)$ of projections in $\cal R$ such that $\sum_{a=1}^ne^x_a=1$ for all $x\in [k]$.  In order to prove their theorem, Kim, Paulsen, and Schaufhauser observed the following fact:

\begin{lem}
For any $k$ and $n$, $Z(k,n)$ is a definable subset of $\cal R$.
\end{lem}

Of course, Kim, Paulsen, and Schaufhauser stated their result in operator-algebraic terms, representing a beautiful confluence of operator-algebraic and model-theoretic ideas!  Consequently, we have:

\begin{cor}
For any nonlocal game $\frak G$ with $k$ questions and $n$ answers, we have
$$\pval^*(\frak G)=\left(\sup_{z\in Z(k,n)}\varphi(z)\right)^{\R}.$$
\end{cor}

An analysis of how formulae with quantifiers over a definable set can be equivalently expressed using ordinary formulae shows that the sentence appearing above is effectively equivalent to a universal sentence, uniform in the description of the game.

The upshot of all of this is that $\operatorname{MIP}^*=\operatorname{RE}$ implies that $\Th_\forall(\R)$ cannot be effectively approximated from above, that is, there is no algorithm which, upon input a universal sentence $\sigma$, returns a computable sequence of upper bounds converging to $\sigma^\R$.  However, a few years prior to the appearance of $\operatorname{MIP}^*=\operatorname{RE}$, the author and Bradd Hart observed the following in \cite{GH}:

\begin{thm}[G. and Hart]
If CEP has a positive answer, then $\Th_\forall(\R)$ \emph{is} effectively approximable from above.
\end{thm}

The proof of the previous theorem proceeds as follows.  If $\sigma$ is a universal sentence in the language of tracial von Neumann algebras, then CEP implies that the supremum of $\sigma^M$ as $M$ ranges over all II$_1$ factors is simply $\sigma^\R$.  On the other hand, the Completeness theorem for continuous logic \cite{BYP} implies that this supremum is equivalently expressed as the infimum of all dyadic rational numbers $r$ such that $T_{II_1}\vdash \sigma\dminus r$, where $T_{II_1}$ is the theory of II$_1$ factors.  Since $T_{II_1}$ is an effectively enumerable theory, the collection of such $r$'s is itself effectively enumerable, implying the desired result.

Putting all of the pieces together, we arrive at the model-theoretic explanation for why a negative solution to CEP follows from $\operatorname{MIP}^*=\operatorname{RE}$.

\subsection{Bells and whistles}

The above analysis leads to a \emph{G\"odelian}-style refutation of the CEP.  To state this precisely, we introduce the following notion.

\begin{defn}
Let $M$ be an structure in some (continuous) language $L$.  By the \textbf{$M$EP} we mean the statement that there is an effectively enumerable subset $T\subseteq \Th(M)$ such that all models of $T$ embed into an ultrapower of $M$.
\end{defn}

The analysis of the previous section shows the following:

\begin{cor}
The $\R$EP has a negative solution.
\end{cor}

Since $\R$ embeds into every II$_1$ factor, we could drop the requirement in the above definition that $T$ be contained in $\Th(M)$ and instead merely ask that $T$ extend the theory of II$_1$ factors, obtaining the following corollary:

\begin{cor}\label{godel}
There is no effectively enumerable theory $T$ extending the theory of II$_1$ factors all of whose models embed into an ultrapower of $\R$.
\end{cor}

In Subsection \ref{mono} below, we will see how to derive the failure of the $A$EP for certain \cstar-algebras $A$ from the above arguments.

A consequence of Corollary \ref{godel} above is that one can infer the existence of infinitely many universal theories of II$_1$ factors, something the ``standard'' refutation of CEP does not seem to imply: 

\begin{cor}
There is a sequence $M_1,M_2,\ldots,$ of separable II$_1$ factors, none of which embed into an ultrapower of $\cal R$, satisfying, for all $i<j$, that $M_i$ does not embed into an ultrapower of $M_j$.  In particular, there are infinitely many universal theories of II$_1$ factors.
\end{cor}

To construct the above sequence, let $M_1$ be any separable II$_1$ factor that does not embed into an ultrapower of $\R$ and let $\sigma_1$ be a universal sentence such that $\sigma_1^{\R}=0$ but $\sigma_1^{M_1}>0$.  Let $r_1\in (0,\sigma_1^{M_1})$ be a rational number and let $T_1:=T_{II_1}\cup \{\sigma_1\dminus r_1\}$.  Since $T_1$ is effectively enumerable, there is a separable model $M_2$ of $T_1$ that does not embed into an ultrapower of $\R$.  Since $\sigma_1^{M_1}>\sigma^{M_2}_1$, we have that $M_1$ does not embed into an ultrapower of $M_2$.  Since $M_2$ does not embed into an ultrapower of $\R$, there is a universal sentence $\sigma_2$ such that $\sigma_2^{\R}=0$ but $\sigma_2^{M_2}>0$.  Let $r_2\in (0,\sigma_2^{M_2})$ be rational and let $T_2:=T_1\cup\{\sigma_2\dminus r_2\}$.  Once again, $T_2$ is effectively enumerable, so there is a separable model $M_3$ of $T_2$ that does not embed into an ultrapower of $\R$.  It is clear that neither $M_1$ nor $M_2$ embed into an ultrapower of $M_3$.  One constructs the remainder of the sequence analogously.

We believe the following question should have a positive answer:

\begin{question}
Do there exist continuum many universal theories of II$_1$ factors?
\end{question}

Another application of Corollary \ref{godel} above is the following, which also appears not to follow from the ``standard'' refutation of the CEP:

\begin{cor}\label{nongamma}
There is a II$_1$ factor without property Gamma that does not embed into an ultrapower of $\R$.
\end{cor}

To see this, simply let $\sigma$ be one of the sentences in the axiomatization of property Gamma for which $\sigma^{L(\bb F_2)}>0$ and let $T:=T_{II_1}\cup \{\sigma\dminus r\}$ for some rational number $r\in (0,\sigma^{L(\bb F_2)})$.  By Corollary \ref{godel} above, there is a model of $T$ that does not embed into an ultrapower of $\R$, which is thus the desired factor.  (Corollary \ref{nongamma} also follows from a very recent result of Chifan, Drimbe, and Ioana \cite{CDI}, who prove that every II$_1$ factor embeds into a II$_1$ factor with property (T).)

The following question is surprisingly more difficult:

\begin{question}
Are there two separable II$_1$ factors $M$ and $N$, each of which embed into an ultrapower of $\R$, neither of which have property Gamma, for which $M\not\equiv N$?
\end{question}

\subsection{The existence of the enforceable II$_1$ factor}

In this subsection, we mention a model-theoretic variant of the CEP that is still open and, in this author's opinion, is one of the more interesting open problems in the model theory of operator algebras.  

As discussed in the article by the author and Hart in this volume, $\R$ is an e.c. model of its universal theory and the CEP is thus equivalent to the assertion that $\R$ is an e.c. II$_1$ factor.  The ideas in this subsection elaborate further on this observation.  We first note the following:

\begin{lem}
Being hyperfinite is a $\forall\bigvee\exists$-axiomatizable property of II$_1$ factors.
\end{lem}

The proof of this lemma has not appeared explicitly in the literature but is similar to the proof of the main results in \cite{AF}.

Armed with this and Proposition \ref{force} above, we arrive at the following:

\begin{cor}
$\R$ is the enforceable model of its universal theory.
\end{cor}

In particular, if CEP were to hold, then $\R$ would be the enforceable II$_1$ factor.  (Of course, the converse is also true, but moot at this point.)  Nevertheless, the following question is open and tantalizing:

\begin{question}
Does the enforceable II$_1$ factor exist?
\end{question}

How likely is it that the enforceable II$_1$ factor exists?  That is of course difficult to say.  If it did exist, then it would ``rival'' $\R$ for being ``the most important II$_1$ factor'' for it would be generic from the model-theoretic point of view.  In \cite{resembles}, some properties of the enforceable II$_1$ factor were established (of course presuming its existence).

It is worth pointing out one notable case when the enforceable object does not exist, namely for the (classical) theory of groups.  (This seems implicit in Hodges' book \cite{hodges} but is written down explicitly in the article \cite{sriyash} by Kunnawalkam Elayavalli, Lodha, and the author.)  However, the ingredients involved in the proof are a blend of recursion-theoretic and combinatorial group-theoretic tools that seem to be currently unavailable to us in the II$_1$ factor setting.

Another interesting variant of the CEP is the following:

\begin{question}
Is the property of being isomorphic to a group von Neumann algebra an enforceable property?
\end{question}

\section{The Kirchberg Embedding Problem}\label{kepsection}

In this and the next section, we consider \cstar-algebra versions of the CEP.  For simplicity, \textbf{henceforth all \cstar-algebras will be assumed to be unital}.  (Much of what is said below can be adapted to the not necessarily unital situation, but this assumption simplifies the exposition.)

If one phrases the CEP as the statement that every tracial von Neumann algebra embeds into the tracial ultrapower of an \emph{injective} II$_1$ factor, then a natural \cstar-algebra analog of CEP would be to ask whether or not every \cstar-algebra embeds into an ultrapower of a nuclear \cstar-algebra (for a \cstar-algebra $A$ is nuclear if and only if its enveloping von Neumann algebra $A^{**}$ is injective).  By Kirchberg's celebrated theorem \cite{KP}, every separable nuclear \cstar-algebra embeds into the Cuntz algebra $\cal O_2$ (see Szab\'o's article in this volume for the definition of $\cal O_2$), whence it is equivalent to ask whether or not $\cal O_2$ is a locally universal \cstar-algebra.  We refer to this problem as the \textbf{Kirchberg embedding problem} (KEP).  (We attribute this problem to Kirchberg as we first learned of this problem from Ilijas Farah, who in turn first learned of this problem during a discussion with Kirchberg in 2007.  The first mention of this problem in the literature appears to be in the author's article \cite{GS} with Sinclair.)

At the moment of the writing of this article, the KEP remains an open problem.  In this section, we mention the connection between it and the model theory of \cstar-algebras in a way that parallels the situation with the CEP.  Recall that the CEP is equivalent to the statement that $\R$ is an e.c. tracial von Neumann algebra.  The analogous statement for the KEP holds:

\begin{prop}\label{kepec}
The KEP has a positive solution if and only if $\O_2$ is an e.c. \cstar-algebra.
\end{prop}

The proof is analogous to the proof of the same statement for tracial von Neumann algebras, using the fact that any two embeddings of $\O_2$ into its ultrapower are unitarily conjugate (a consequence of the fact that $\O_2$ is a \textbf{strongly self-absorbing} \cstar-algebra) and that the theory of \cstar-algebras has the joint embedding property.

The following variation on the preceding proposition is also of interest:

\begin{prop}\label{only}
$\O_2$ is the only possible separable \cstar-algebra that is nuclear and e.c.  Consequently, a positive solution to the KEP is equivalent to the statement that there is a \cstar-algebra that is both nuclear and e.c.
\end{prop}

The proof of this proposition is quite interesting.  Indeed, suppose that $A$ is a separable \cstar-algebra that is both nuclear and e.c.  A consequence of being e.c. is that $A$ is simple (see \cite{GS}).  By another fundamental result of Kirchberg \cite{KP}, the fact that $A$ is simple, separable and nuclear implies that $A\otimes \O_2\cong \O_2$.  However, a consequence of $A$ being e.c. is that $A$ is ``$\O_2$-stable,'' that is, $A\otimes \O_2\cong A$; this follows from the fact that being $\O_2$-stable is an $\forall\exists$-axiomatizable property of \cstar-algebras (see \cite{munster}) together with the fact that every \cstar-algebra embeds into an $\O_2$-stable algebra (namely by tensoring the algebra with $\O_2$ itself and using the fact, due to Cuntz, that $\O_2\otimes \O_2\cong\O_2$).  It follows that $A\cong \cal O_2$, as desired.  An alternative, slightly more elementary proof, can be found in \cite[Remark 21]{international}.

Another similarity with the CEP concerns enforceability.  First, we will need the following result:

\begin{prop}
Being nuclear is a $\forall\bigvee\exists$-axiomatizable property of \cstar-algebras. 
\end{prop}

Two proofs for the preceding proposition are offered in \cite{munster}, one ``soft'' and model-theoretic, the other ``concrete,'' writing down specific axioms for nuclearity.  Let us sketch the former argument.  Given a \cstar-algebra $A$ and $k,n\in \bb N$, define the predicate $\nuc_{k,n}^A:A_1^k\to \bb R$ by $\nuc_{k,n}^A(\vec a)=\inf_{\phi,\psi}\|(\psi\circ \phi)(\vec a)-\vec a\|$, where $\phi$ ranges over all ucp maps $A\to M_n(\bb C)$ and $\psi$ ranges over all ucp maps $M_n(\bb C)\to A$.  An argument using the Beth definability theorem shows that the predicates $\nuc_{k,n}$ are actually existentially definable relative to the theory of \cstar-algebras, that is, there are existential formulae (really, uniform limits of existential formulae) $\Phi_{k,n}(\vec x)$ such that, for every \cstar-algebra $A$, every $k,n\in \bb N$, and every $\vec a\in A^k_1$, we have $\nuc_{k,n}^A(\vec a)=\Phi_{k,n}^A(\vec a)$.  It remains to note that a \cstar-algebra $A$ is nuclear if and only if, for every $k\in \bb N$, we have $\left(\sup_{\vec x}\bigvee_{n\in \bb N}\Phi_{k,n}(\vec x)\right)^A=0$.

As with the CEP, we arrive at the following formulation of the KEP, whose proof uses everything we have discussed thus far:

\begin{thm}[G. \cite{enforceable}]
The following are equivalent:
\begin{enumerate}
    \item The KEP has a positive solution.
    \item The property of being nuclear is enforceable.
    \item $\O_2$ is the enforceable \cstar-algebra.
\end{enumerate}
\end{thm}

The previous theorem yielocus an interesting local, finitary reformulation of the KEP first identified by Sinclair and the author in \cite{GS}.  First, we say a condition (relative to the theory of \cstar-algebras) $p(\vec x)$, where $\vec x$ is a $k$-tuple, has \textbf{good nuclear witnesses} if, for every $\epsilon>0$, there is a \cstar-algebra $A$, $\vec a \in A^k$, and an $n\in \bb N$ such that $\vec a$ satisfies the condition $p(\vec x)$ and for which $\nuc_{k,n}^A(\vec a)<\epsilon$.  The previous theorem then yields the following corollary:

\begin{cor}
The KEP has a positive solution if and only if every condition has good nuclear witnesses.
\end{cor}

The import of the previous corollary is a (seemingly) significant weakening of the demand that every condition be satisfied in a nuclear \cstar-algebra (equivalently, satisfied in $\O_2$), for one only asks that the witness admit a good ucp factorization through a matrix algebra, and, moreover, the witness and the dimension of the matrix algebra can vary as the level of approximation varies.

\section{The MF problem and the quasidiagonality problem}\label{MFsection}

\subsection{A negative solution to the MF problem}

Recall that the CEP is also equivalent to the statement that every separable tracial von Neumann algebra embeds into $\prod_\u M_n(\mathbb C)$, a \emph{tracial} ultraproduct of matrix algebras with respect to a nonprincipal ultrafilter on $\bb N$.  It is thus natural to formulate a \cstar-algebra version of CEP by asking that every separable \cstar-algebra embed into a \emph{\cstar-algebra} ultraproduct $\prod_\u M_n(\bb C)$ of matrix algebras with respect to a nonprincipal ultrafilter on $\bb N$.  Using the same notation for both tracial von Neumann algebra ultraproducts and \cstar-algebra ultraproducts is potentially dangerous (and some authors even use different notations for the two ultraproducts); to prevent confusion, \textbf{in the remainder of this section, unless explicitly stated otherwise, all ultraproducts will be \cstar-algebra ultraproducts}.

There is an immediate obstruction to the statement ``every separable \cstar-algebra embeds into $\prod_\u M_n(\bb C)$'' from being true, namely the \cstar-algebra $\prod_\u M_n(\bb C)$ is stably finite, as is any subalgebra.  Thus, we may modify the problem as follows:

\begin{defn}
The \textbf{MF problem} is the problem of whether or not every separable stably finite \cstar-algebra embeds into $\prod_\u M_n(\bb C)$.
\end{defn}

The terminology MF comes from the fact that a separable \cstar-algebra is called \textbf{matricially finite} (or MF) if it embeds into $\prod_\u M_n(\bb C)$.  Consequently, the MF problem asks if the notions of stably finite and MF coincide for separable \cstar-algebras.

As in the case of the CEP, the MF problem can be reformulated in terms of ultrapowers of a single object.  Indeed, the MF problem is equivalent to the problem of whether every separable stably finite \cstar-algebra embeds into a nonprincipal ultrapower $\cal Q^\u$ of the \textbf{universal UHF algebra} $\cal Q$ (see \cite[Lemma 4.4.1]{munster}).

An immediate consequence of the negative solution of the CEP is that the MF problem also has a negative solution:

\begin{cor}
The MF problem has a negative solution.
\end{cor}

To prove the previous corollary, suppose that $M$ is a II$_1$ factor that does not embed into $\R^\u$ (here we mean the tracial von Neumann algebra ultrapower).  We claim then that $M$ does not embed (as a C*-algebra) into a nonprincipal C*-ultrapower of $\cal Q$; by considering a separable elementary subalgebra of $M$ (in the language of \cstar-algebras), we obtain the desired counterexample to the MF problem.  Suppose, towards a contradiction, that $i:M\hookrightarrow \cal Q^\u$ is an embedding.  
%(To be clear, $\cal Q^\u$ is the C*-algebra ultrapower of $\cal Q$.).
Recall that $\cal Q$ has a unique trace $\tau$ and the von Neumann algebra generated by $\cal Q$ with respect to the GNS represenation corresponding to $\tau$ is $\cal R$.  Let $\pi:\mathcal Q^\u \to \mathcal R^\u$ denote the composition of the quotient map $\cal Q^\u\to \cal Q^\u/I$, where $$I=\{(x_n)^\bullet \in \cal Q^\u \ : \ \lim_\u\|x_n\|_2=0\}$$ is the trace ideal, with the natural inclusion $\cal Q^\u/I\hookrightarrow \R^\u$ obtained by viewing operator norm bounded balls in $\cal Q$ as $\|\cdot\|_2$-dense subsets of the corresponding balls in $\cal R$.  Since $M$ has a unique trace, which is faithful, we get that the composition  $\pi\circ i:M\to \mathcal R^{\u}$ is a trace-preserving *-homomorphism, a contradiction.

The following question seems wide open:

\begin{question}
Does $\R$ embed, as a \cstar-algebra, into $\cal Q^\u$?
\end{question}

The negative solution to the MF problem also has a G\"odelian-style refutation, that is, the $\cal Q$EP has a negative solution as well; we postpone the discussion of this fact until Subsection \ref{mono} below.

Unlike most of the embedding problems discussed in this paper, it is not even clear that there \emph{ought} to be a locally universal object for the class of stably finite \cstar-algebras as the following question appears to be open:

\begin{question}
Does the class of stably finite \cstar-algebras have the JEP?
\end{question}

A natural guess would be that the minimal tensor product of two stably finite \cstar-algebras would once again be stably finite.  However, the validity of this statement is far from clear.  In fact, the question of whether or not the minimal tensor product of two \emph{simple} stably finite \cstar-algebras is once again stably finite is equivalent to a well-known open problem, namely whether or not every stably finite \cstar-algebra admits a trace (see \cite{international}).  It is worth mentioning that the class of \cstar-algebras admitting a trace does have JEP (the tensor product trace on the minimal tensor product witnesses this) and consequently any e.c. object for this class (which exists since the class is inductive) is locally universal.

\subsection{The quasidiagonality problem}

Connes' original motivation for considering the question of which tracial von Neumann algebras embed into ultrapowers of $\R$ came from his striking result proving that injective II$_1$ factors were hyperfinite, thus completing the classification of injective II$_1$ factors \cite{connes}.  A crucial ingredient in his proof was that injective factors did indeed admit embeddings into ultrapowers of $\R$.

In \cstar-algebra theory, the analogous problem would be trying to classify simple, nuclear \cstar-algebras, where simple is the analog of being a factor and (as already mentioned) nuclear is the analog of being injective.  In trying to mimic Connes' approach in classifying simple, nuclear \cstar-algebras, it thus becomes natural to try to prove that they admit embeddings into $\cal Q^\u$.  As stated in the previous subsection, an immediate obstruction to proving such a result is that the \cstar-algebra in question must be stably finite.  Thus, one is naturally led to:

\begin{defn}[Quasidiagonality problem-simple version]
Does every simple, stably finite, nuclear \cstar-algebra embed into $\cal Q^\u$?
\end{defn}

A word about the nomenclature in the previous definition is in order.  A \cstar-algebra $A$ is said to be \textbf{quasidiagonal} if there is an embedding $A\hookrightarrow \prod_\u M_n(\bb C)$ that admits a ucp lift $A\to \prod_{n\in \bb N}M_n(\bb C)$.  Thus, quasidiagonal \cstar-algebras form a special subclass of the class of MF-algebras.  However, by the Choi-Effros lifting theorem, if a \emph{nuclear} \cstar-algebra is MF, then the aforementioned ucp lift automatically exists, whence there is no difference in the two notions.  Halmos defined what it meant for a set of bounded operators to be quasidiagonal and then a \cstar-algebra is called quasidiagonal if it admits a concrete representation for which the operators in the image of the representation form a quasidiagonal set.  Voiculescu then proved that this definition of quasidiagonal \cstar-algebra agrees with the one given at the beginning of this paragraph (see \cite{voiculescu}).

As with the MF problem, it is not evident that the class of simple, stably finite nuclear \cstar-algebras should have a locally universal object for this class is not known to have JEP.

% We should point out that there is a locally universal object for the category of simple, stably finite, nuclear \cstar-algebras:

% \begin{prop}
% There are separable e.c. objects in the class of simple, stably finite, nuclear \cstar-algebras.  Any such object is locally universal for this class.
% \end{prop}

% The first part of the previous proposition follows simply from the fact that the class of separable, simple, stably finite, nuclear \cstar-algebras is closed under countable directed unions.  To see the second assertion, it suffices to...ACTUALLY IS THIS TRUE?  MAYBE MENTION THE TRACIAL VERSION, WHERE THEY DO EXIST.

Amazingly enough, this modified version of the MF problem has almost been shown to be true, modulo one technical assumption:

\begin{thm}[Tikuisis, Winter, White \cite{TWW}]\label{tww}
Every simple, stably finite, nuclear \cstar-algebra \emph{satisfying the UCT} is quasidiagonal.
\end{thm}

Here, the UCT is short for the \textbf{Universal Coefficient Theorem}.  Assuming that a \cstar-algebra satisfies the UCT is a technical K-theoretic assumption on the algebra.  (See \cite{rordam} for more information on the UCT.)  One of the major open questions in \cstar-algebra theory is:  

\begin{question}[UCT problem]
Do all separable nuclear \cstar-algebras satisfy the UCT? 
\end{question} 

Of all of the adjectives appearing in the statement of Theorem \ref{tww}, all but the UCT have been shown to have model-theoretic meaning:  being stably finite and MF are universally axiomatizable properties whilst being simple and being nuclear are $\forall\bigvee \exists$-axiomatizable properties.  (It turns out that quasidiagonality in general is also $\forall\bigvee\exists$-axiomatizable; see \cite[Section 5.13]{munster}.)  

It is interesting to ask:

\begin{question}
Is satisfying the UCT a $\forall\bigvee\exists$-axiomatizable property of separable nuclear \cstar-algebras?.
\end{question}   

The previous question notwithstanding, Barlak and Szab\'o \cite{BS} proved the following interesting fact:

\begin{thm}\label{bs}
If $A$ is an e.c. subalgebra of $B$ and $B$ is a nuclear \cstar-algebra satisfying the UCT, then so does $A$.
\end{thm}

This theorem allows one to deduce the truth of the simple version of the quasidiagonality problem from a weakening of the UCT problem:

\begin{thm}[G. and Sinclair \cite{international}]\label{qdsimple}
Suppose that every simple, stably finite, nuclear \cstar-algebra embeds into a simple, stably finite, nuclear \cstar-algebra satisfying the UCT.  Then the simple version of the quasidiagonality problem is true.
\end{thm}

Indeed, suppose that $A$ is an e.c. simple, stably finite, nuclear \cstar-algebra (which exists since this class is $\forall\bigvee\exists$-axiomatizable).  By the assumption of the theorem and Theorem \ref{bs} above, it follows that $A$ itself satisfies the UCT, whence $A$ is quasidiagonal by Theorem \ref{tww} above.  Now if $B$ is any separable, simple, stably finite, nuclear \cstar-algebra, then $A\otimes B$ is also stably finite; this uses the quasidiagonality of $A$ (see \cite[Lemma 4]{international}).  Since $A$ is e.c., we have that $A\otimes B$, and thus $B$, embeds into an ultrapower of $A$, whence $A$ is locally universal for the class of simple, stably finite, nuclear \cstar-algebras.  Since $A$ is quasidiagonal, the result follows.

One can remove the simplicty assumption in the quasidiagonality problem, arriving at:

\begin{defn}[Quasidiagonality problem-general version]
Does every stably finite, nuclear \cstar-algebra embed into $\cal Q^\u$?
\end{defn}

The quasidiagonality problem has several model-theoretic equivalents:

\begin{thm}
The following are equivalent:
\begin{enumerate}
    \item The quasidiagonality problem has a positive solution.
    \item Being UHF is an enforceable property of stably finite, nuclear \cstar-algebras.
    \item $\cal Q$ is the enforceable stably finite, nuclear \cstar-algebra.
    \item $\cal Q$ is an e.c. stably finite, nuclear \cstar-algebra.
\end{enumerate}
\end{thm}

The implication (1) implies (2) in the previous proposition follows from Proposition \ref{force} and the fact that being UHF is $\forall\bigvee\exists$-axiomatizable (see \cite{AF}).  The implication (3) implies (4) follows from the fact that being e.c. is enforceable while the implication (4) implies (1) proceeds along the lines of Theorem \ref{qdsimple} above.  Finally, to see (2) implies (3), one first notes that the property of being $\cal Q$-stable, that is, that $A\otimes \cal Q\cong A$, is also enforceable.  Indeed, this property is $\forall\exists$-axiomatizable (for the same reason as in the case of $\cal O_2$) and is thus true of any e.c. object in this class as any object in this class is a subalgebra of a $\cal Q$-stable object in the class (by tensoring with $\cal Q$).  It remains to note that $\cal Q$ is the only $\cal Q$-stable UHF algebra.

The following question is the natural stably finite analog of Theorem \ref{only} above; see \cite{international} for partial progress towards its resolution:

\begin{question}
Suppose that $A$ is an e.c. stably finite \cstar-algebra that is also nuclear.  Must we have $A\cong \cal Q$?
\end{question}

At the moment, it is unclear if the general version of the quasidiagonality problem could be deduced from the simple version, even assuming a positive solution to the UCT problem.  However, using model-theoretic forcing again, one can prove the following result:

\begin{thm}[G. and Sinclair \cite{international}]
Suppose the following hold:
\begin{enumerate}
    \item Every stably finite nuclear \cstar-algebra embeds into a stably finite, nuclear \cstar-algebra satisfying the UCT.
    \item There is a simple, stably finite, nuclear \cstar-algebra that is locally universal for the class of stably finite nuclear \cstar-algebras. 
\end{enumerate}
Then the general version of the quasidiagonality problem holds.
\end{thm}

Note that the second item in the hypotheses of the previous theorem is indeed a weakening of the statement of the quasidiagonality problem as $\cal Q$ itself is simple.  The proof of the preceding theorem proceeds similarly as in the proof of Theorem \ref{qdsimple} above.  Indeed, the second condition, Proposition \ref{force}, and the fact that being simple is $\forall\bigvee\exists$-axiomatizable allows one to construct an e.c. stably finite nuclear \cstar-algebra $A$ that is simple.  Moreover, the first condition and Theorem \ref{bs} above allows one to conclude that $A$ satisfies the UCT.  Thus, Theorem \ref{tww} above allows one to conclude that $A$ is quasidiagonal.  It follows that the quasidiagonality problem has a positive answer just as in the conclusion of the proof of Theorem \ref{qdsimple}.

By the negative solution to the MF problem, there is a universal sentence $\sigma$ in the language of \cstar-algebras for which $\sigma^{\cal Q}=0$ and yet $\sigma^A=r>0$ for some stably finite \cstar-algebra $A$.  Let $T$ be the theory of stably finite \cstar-algebras together with the existential condition $r\dminus \sigma=0$.  We adapt the terminology from Section \ref{kepsection} above and say that a condition $p(\vec x)$ (with $\vec x=(x_1,\ldots,x_k)$) relative to the theory $T$ has good nuclear witnesses if, for every $\epsilon>0$, there is $A\models T$ and $\vec a\in A$ satisfying $p$ such that $\nuc_{k,n}^A(\vec a)<\epsilon$.

\begin{thm}
Using the terminology in the previous paragraph, suppose that every condition has good nuclear witnesses.  Then the quasidiagonality problem has a negative solution.
\end{thm}

Indeed, the assumption that every condition has good nuclear witnesses allows one to construct a model of $T$ that is nuclear; being a model of $T$, the algebra is also stably finite but not embeddable in an ultrapower of $\cal Q$ (that is, not quasidiagonal).  Considering the contrapositive of this theorem, if the quasidiagonality problem has a positive solution, then whenever one has a counterexample to the MF problem as in the paragraph above, then there must be some condition $p(\vec x)$ relative to the associated theory $T$ which does not have good nuclear witnesses, meaning that there is some some $\epsilon>0$ such that, in every model $A$ of $T$, every $k$-tuple $\vec a$ from $A$ satisfying $p$ must satisfy $\nuc_{k,n}^A(\vec a)\geq \epsilon$ for all $n\in \bb N$.  

\subsection{Monotracial \cstar-algebras and the Jiang-Su Embedding Problem}\label{mono}

One can adapt the techniques used to show that the $\R$EP fails to show that the $A$EP fails for a large class of \cstar-algebras.  We will be concerned with \textbf{monotracial} \cstar-algebras, that is, \cstar-algebras which admit a unique tracial state.  For example, the universal UHF algebra $\cal Q$ is monotracial.  Suppose that $A$ is a monotracial \cstar-algebra whose unique trace is $\tau_A$.  Let $N$ denote the von Neumann algebra generated by $A$ via the GNS representation of $A$ associated to $\tau_A$.  Then $N$ is a tracial von Neumann algebra when equipped with the extension $\tau_N$ of the original trace $\tau_A$ to $N$.  Since $A$ has a unique trace, it follows that $\tau_N$ is the unique trace on $N$, whence $N$ is in fact a tracial factor.  If we further assume that $A$ is infinite-dimensional, then we can conclude that $N$ is a II$_1$ factor.

Now suppose that $\sigma$ is a universal sentence in the language of tracial von Neumann algebras.  We can also view $\sigma$ as a sentence in the language of tracial \cstar-algebras (one that simply does not refer to the operator norm in any way) and can thus compare the values $\sigma^{(A,\tau_A)}$ and $\sigma^{(N,\tau_N)}$.  Since $A$ is $\|\cdot\|_2$-dense in any operator norm bounded subset of $N$, it is clear that $\sigma^{(N,\tau_N)}\leq \sigma^{(A,\tau_A)}$.  On the other hand, if $A$ is further assumed to be simple, then $A$ embeds into $N$ and thus, $\sigma^{(A,\tau_A)}=\sigma^{(N,\tau_N)}$.

Summarizing thus far:  if $A$ is a unital, simple, infinite-dimensional, monotracial \cstar-algebra whose weak closure in the GNS representation we denote by $N$, then for any sentence $\sigma$ in the language of tracial von Neumann algebras, we have that $\sigma^{(A,\tau_A)}=\sigma^{(N,\tau_N)}$.  If we further assume that $N$ embeds in an ultrapower of $\R$, then this common value equals $\sigma^{(\R,\tau_R)}$.

It is thus tempting to try to conclude that the $A$EP must fail for any \cstar-algebra $A$ satisfying the conditions appearing in the previous paragraph.  Indeed, if, towards a contradiction, there was an effectively enumerable subset $T\subseteq \Th(A)$, all of whose models embed into an ultrapower of $A$, then by adding to these axioms the (effectively enumerable) axioms for tracial \cstar-algebras, one might hope that by running proofs from this new theory $T'$, one might be able to obtain effective upper bounds for $\Th_\forall(\R)$ and thus contradict $\operatorname{MIP}^*=\operatorname{RE}$ as before.  The issue with this is that if $(M,\tau)\models T'$, one is only guaranteed that $M$ embeds into $A^\u$ as a \cstar-algebra, that is, the embedding need not preserve the trace $\tau$ on $M$.  Consequently, we would not know that $\sup\{\sigma^{(M,\tau)} \ : \ (M,\tau)\models T'\}$ coincides with $\sigma^{(A,\tau_A)}=\sigma^{(\R,\tau_R)}$ and thus our usual Completeness Theorem argument need not go through.  However, if $M$ were itself monotracial, then the above embedding would be guaranteed to be trace-preserving and the above argument would work.  Unfortunately, being monotracial is not an axiomatizable property of \cstar-algebras.  That being said, there is a an axiomatizable property of \cstar-algebras known as the \textbf{uniform Dixmier property} which implies being monotracial.  (To be fair, the uniform Dixmier property itself is not axiomatizable.  Instead, there are quantitative versions known as the $(m,\gamma)$-uniform Dixmier property for some parameters $m\in \bb N$ and $\gamma\in (0,1)$, which are each axiomatizable; having the uniform Dixmier property means having the $(m,\gamma)$-Dixmier property for some choice of parameters $m$ and $\gamma$.  See \cite{art} for details.)

In summary, we have:

\begin{thm}[G. and Hart \cite{undec}]\label{AEP}
Suppose that $A$ is a unital, infinite-dimensional, simple, \cstar-algebra $A$ with the uniform Dixmier property whose associated GNS von Neumann algebra $N$ embeds in an ultrapower of $\R$.  Then the $A$EP has a negative solution.
\end{thm}

There are many examples of \cstar-algebras satisfying the hypotheses appearing in the previous theorem.  In particular, it follows from \cite{hz} and \cite[Corollary 3.11]{art} that $\cal Q$ satisfies all of the above hypotheses, leading to the aforementioned strengthed refutation of the MF problem:

\begin{cor}
The $\cal Q$EP has a negative solution.
\end{cor}

We can use Theorem \ref{AEP} to prove a new, purely operator algebra-theoretic result, which refutes another natural ultarpower embedding problem.

One of the most important algebras in modern \cstar-algebra classification theory is the \textbf{Jiang-Su} algebra $\cal Z$.  (See Vignati's article in this volume for more information on $\cal Z$.)  It follows from the works in \cite{hz} and \cite[Remark 3.18 and Corollary 3.22]{art} that $\cal Z$ also satisfies hypotheses of Theorem \ref{AEP}, whence we have:

\begin{cor}\label{ZEP}
The $\cal Z$EP has a negative solution.
\end{cor}

One of the defining features of $\cal Z$ is that it is \textbf{stably projectionless}, meaning that for any $n\in \bb N$ and any projection $p\in M_n(\cal Z)$, there is a projection $q\in M_n(\bb C)$ unitarily conjugate to $p$.  Being stably projectionless is axiomatizable by the following (effective) list of axioms, one for each $n\in \bb N$:  $$\sup_{p\in M_n( A)}\inf_{q\in M_n(\bb C)}\inf_{u\in U(M_n(A))}d(upu^*,q)=0.$$  A couple of words are in order about this axiomatization.  First, it is know that the matrix amplifications $M_n(A)$ belong to the imaginary sorts of the theory of \cstar-algebras, whence the first and third quantifiers are not problematic.  Similarly, being (locally) compact, adding the matrix algebras $M_n(\bb C)$ to the theory of \cstar-algebras is also harmless.  Finally, the above axioms only seem to say that every projection in $M_n(A)$ is approximately unitarily equivalent to a projection in $M_n(\bb C)$; however, it is well-known that two projections that are sufficiently close are actually unitarily conjugate, whence the axioms do indeed express that an algebra is stably projectionless.  

Combining Corollary \ref{ZEP} with the discussion in the previous paragraph yields the following stably projectionless analog of the negative solution to the MF problem, which currently has no purely operator-algebraic proof: 

\begin{cor}[G. and Hart \cite{undec}]
There is a stably projectionless \cstar-algebra that does not embed into an ultrapower of $\cal Z$.
\end{cor}
\section{Abelian \cstar-algebras}\label{abeliansection}

In this section, we move on from the more difficult embedding problems in the previous sections and instead consider the case of abelian \cstar-algebras.  Not surprisingly, the ensuing discussion becomes topological in nature.  

As in the previous two sections, all \cstar-algebras in this section are assumed to be unital.

\subsection{Preliminaries on ultracoproducts of compact Hausdorff spaces} 
Some of the model theory of abelian \cstar-algebras to be discussed below follows immediately from classical model-theoretic facts about Boolean algebras together with a categorical understanding of the relevant ultraproduct constructions.  The ideas presented here are an elaboration of those in \cite[Section 5]{EV}.  

To begin, let $\zdcomp$ and $\bool$ denote the categories of zero-dimensional compact Hausdorff spaces and Boolean algebras respectively. Consider the \textbf{Stone functor} $\zdcomp\to \bool$ given by sending the zero-dimensional compact Hausdorff space $X$ to the Boolean algebra $Cl(X)$ of clopen subsets of $X$.  This functor is contravariant and is a duality of categories whose inverse is given by the functor taking a Boolean algebra $\bb B$ to its spectrum, that is, the set of ultrafilters on $\bb B$, or, equivalently, the set of Boolean algebra homomorphisms $\bb B\to \{0,1\}$.  

Letting $\comp$ and $\ab$ denote the categories of compact Hausdorff spaces and unital abelian \cstar-algebras, then we may also consider the \textbf{Gelfand functor} $\comp\to\ab$ given by sending the compact Hausdorff space $X$ to the unital abelian \cstar-algebra $C(X)$ of complex-valued continuous functions on $X$.  Like the Stone functor, the Gelfand functor is contravariant and is a duality of categories, this time the inverse given by the functor taking a unital abelian \cstar-algebra $A$ to its spectrum $\Sigma(A)$ consisting of all unital $*$-homomorphisms $A\to \bb C$.  

In a sense, the Gelfand functor is an ``extension'' of the Stone functor.  More precisely, recall first that a unital \cstar-algebra is called \textbf{real rank zero} if the set of invertible self-adjoint elements is dense in the set of self-adjoint elements.  A unital abelian \cstar-algebra $C(X)$ is real rank zero if and only if $X$ is zero-dimensional.  Let $\rrz$ denotes the category of real-rank zero unital abelian \cstar-algebras.  Then the covariant functor $\rrz\to \bool$ given by composing the inverse of the Gelfand functor (restricted to $\rrz$) with the Stone functor is an equivalence of categories.  In this case, for $X\in \zdcomp$, $C(X)$ is the closed linear span of its space of projections $P(X)$, which in turn is naturally isomorphic to $Cl(X)$.  We refer to this equivalence of categories as the ``forgetful functor'' as it forgets the \cstar-algebra structure and only remembers the Boolean algebra structure on the set of projections.  

Next recall that one can present the \cstar-algebra ultraproduct construction in purely categorical language.  Indeed, suppose that $(A_i)_{i\in I}$ is a family of \cstar-algebras and $\u$ is an ultrafilter on $I$.  For each $J\in \u$, let $A_J:=\prod_{j\in J}A_j$ denote the direct product and note that the family $(A_J,\pi_{JK})$ forms a directed family, where $\pi_{JK}:A_K\to A_J$ is the canonical projection map when $J,K\in \u$ are such that $J\subseteq K$.  There is then a natural isomorphism $\prod_\u A_i\cong \varinjlim A_J$.  This fact is actually completely general and holds for ultraproducts of $L$-structures for any (classical or continuous) language $L$.  In particular, the same observation holds verbatim for Boolean algebras and their ultraproducts.  (See \cite[Chapter 6, Section 10]{ultrabook} for more details.)

Now suppose that $(X_i)_{i\in I}$ is a family of compact Hausdorff spaces and $\u$ is an ultrafilter on $I$.  Since $\prod_\u C(X_i)$ is once again a unital abelian \cstar-algebra, it makes sense to consider the compact Hausdorff space $\Sigma(\prod_\u C(X_i))$.  One can give a purely topological description of $\Sigma(\prod_\u C(X_i))$.  Towards this end, for $J\in \u$, set $X_J:=\coprod_{i\in J} X_i=\beta(\oplus_{i\in J}X_i)$, the Stone-Cech compactification of the direct sum of the $X_i$'s, which is the coproduct construction for the category $\comp$.  Since the Gelfand functor is a duality of categories, it follows that $\prod_{i\in J}C(X_i)\cong C(X_J)$. Applying the inverse of the Gelfand functor to the isomorphism $\prod_\u C(X_i)\cong \varinjlim C(X_J)$ yields the isomorphism $\Sigma(\prod_\u C(X_i))\cong \varprojlim X_J$.  The compact Hausdorff space $\varprojlim X_J$ is called the \textbf{ultracoproduct} of the family $(X_i)_{i\in I}$ with respect to $\u$, denoted $\coprod_\u X_i$.  If each $X_i=X$, we speak of the \textbf{ultracopower} of $X$ with respect to $\u$, denoted $X_\u$.

Suppose now that each $X_i$ in the previous paragraph is also assumed to be zero-dimensional.  One can then apply the forgetful functor to the isomorphisms $$C\left(\coprod_\u X_i\right)\cong \prod_\u C(X_i)\cong \varinjlim C(X_J)$$ to get the isomorphisms $$Cl\left(\coprod_\u X_i\right)\cong \prod_\u Cl(X_i)\cong \varinjlim Cl(X_J).$$  Arguing in a similar fashion, one sees that the forgetful functor sends the diagonal embedding $C(X)\hookrightarrow C(X)^\u$ to the corresponding diagonal embedding $Cl(X)\hookrightarrow Cl(X)^\u$.

\subsection{The model companion of $T_{\ab}$}

The discussion in the previous subsection, together with classical facts about the model theory of Boolean algebras, will allow us to immediately deduce the existence of an $\aleph_0$-categorical model companion for the theory $T_{\ab}$ of unital abelian \cstar-algebras.  We first observe:

\begin{prop}\label{eeforgetful}
Suppose that $X$ and $Y$ are zero-dimensional compact Hausdorff spaces.  Then $C(X)\equiv C(Y)$ if and only if $Cl(X)\equiv Cl(Y)$.
\end{prop}

The shortest proof of the previous proposition appeals to the Keisler-Shelah theorem.  Indeed, if $\u$ and $\cal V$ are ultrafilters, then by the results discussed at the end of the previous subsection, $C(X)^\u\cong C(Y)^{\cal V}$ if and only if $X_\u\cong Y_{\cal V}$ (homeomorphic) if and only if $Cl(X)^\u\cong Cl(Y)^{\cal V}$.

\begin{prop}\label{ecforgetful}
Suppose that $X$ is a compact Hausdorff space.  Then $C(X)$ is an e.c. model of $T_{\ab}$ if and only if (i) $X$ is zero-dimensional, and (ii) $Cl(X)$ is an e.c. Boolean algebra.
\end{prop}

For the forward direction, to prove item (i), we note that being real rank $0$ is $\forall\exists$-axiomatizable (see \cite[Section 3.6.2]{munster}) and every model of $T_{\ab}$ embeds in a real rank zero model of $T_{\ab}$ (since every separable compact metric space is a continuous image of Cantor space).  Item (ii) follows from our analysis in the previous subsection:  if $Cl(X)\subseteq Cl(Y)$, then $C(X)\subseteq C(Y)$ and thus there is an embedding $C(Y)\hookrightarrow C(X)^\u$ that restricts to the diagonal embedding on $C(X)$.  Applying the forgetful functor shows that $Cl(X)$ is e.c. in $Cl(Y)$.  The backwards direction is proven in a similar manner.

We remind the reader that the theory of Boolean algebras has an $\aleph_0$-categorical model completion, namely the theory of atomless Boolean algebras.  In particular, $Cl(2^{\bb N})$ is the unique countable model of this model completion. 

With everything in place, we can now conclude:

\begin{thm}\label{boring}
$\Th(C(2^{\bb N}))$ is $\aleph_0$-categorical and is the model completion of $T_{\ab}$.
\end{thm}

The $\aleph_0$-categoricity of $\Th(C(2^{\bb N}))$ follows from Proposition \ref{eeforgetful} above, the fact that being real-rank zero is elementary, and the $\aleph_0$-categoricity of $\Th(Cl(2^{\bb N}))$. Propositions \ref{eeforgetful} and \ref{ecforgetful}, together with the fact that $\Th(Cl(2^{\bb N}))$ is the model companion of the theory of Boolean algebras, shows that $\Th(C(2^{\bb N}))$ axiomatizes the e.c. models of $T_{\ab}$.  To see that the model companion is in fact a model completion, we simply use the fact that $T_{\ab}$ has the amalgamation property, which follows from the \textbf{fiber product} construction for compact spaces:  if $C(X)$ is embedded in $C(Y)$ and $C(Z)$, then there are surjections $\pi_Y:Y\to X$ and $\pi_Z:Z\to X$ that induce these embeddings.  The corresponding fiber product is the space $W=Y\times_X Z:=\{(y,z)\in Y\times Z \ : \ \pi_Y(y)=\pi_Z(z)\}$.  The natural projection mappings $\theta_Y:W\to Y$ and $\theta_Z:W\to Z$ satisfy $\pi_Y\circ \theta_Y=\pi_Z\circ \theta_Z$.  It follows that the embeddings of $C(Y)$ and $C(Z)$ into $C(W)$ corresponding to $\theta_Y$ and $\theta_Z$ yield the desired amalgamation.  

% \begin{proof}
% If $C(X)$ is separable and $C(X)\equiv C(2^{\bb N})$, then $Cl(X)\equiv Cl(2^{\bb N})$, whence $Cl(X)\cong Cl(2^{\bb N})$ and thus $C(X)\cong C(2^{\bb N})$.

% Since being real rank $0$ is $\forall\exists$-axiomatizable and every abelian \cstar-algebra embeds in a real rank zero abelian \cstar-algebra (since every separable compact metric space is a continuous image of Cantor space), we have that every e.c. abelian \cstar-algebra $C(X)$ is real rank zero.  But then $CL(X)$ is an e.c. Boolean algebra, whence $CL(X)\equiv CL(2^{\bb N})$ and thus $C(X)\equiv C(2^{\bb N})$.  Conversely, if $C(X)\equiv C(2^{\bb N})$, then $C(X)$ has real rank $0$ and $CL(X)\equiv CL(2^{\bb N})$, whence $CL(X)$ is an e.c. Boolean algebra and thus $C(X)$ is an e.c. abelian \cstar-algebra.  Consequently, we have shown that $\Th(C(2^{\bb N}))$ axiomatizes the class of e.c. abelian \cstar-algebras, whence the model companion exists.

% To see that the model companion is in fact a model completion, we simply use the fact that $T_{\ab}$ has the amalgamation property, which follows from the fiber product construction for compact spaces.
% \end{proof}

\subsection{The projectionless case}

While the model theory of the entire class of unital abelian \cstar-algebras is fairly mundane (in the sense that Theorem \ref{boring} above is basically a classical result in disguise), the situation when one considers the subclass of \emph{projectionless} unital abelian \cstar-algebras is far more interesting.  Note that $C(X)$ is projectionless if and only if $X$ is a \emph{connected} compact Hausdorff space, otherwise known as a \textbf{continuum}.  The collection of projectionless unital abelian \cstar-algebras does indeed form an elementary class, which follows either from the observation that the ultracoproduct of a family of continua is once again a continuum (see \cite{gurevic}) or from writing down concrete axioms in the language of \cstar-algebras (see \cite{EGV}).  We let $T_{cont}$ denote the theory of projectionless unital abelian \cstar-algebras and note that this theory is universally axiomatizable.

Unlike the case of arbitrary abelian \cstar-algebras, surprisingly all (nondegenerate) projectionless abelian \cstar-algebras have the same universal theory:

\begin{thm}[K.P. Hart \cite{KP}]\label{sameuniversal}
If $X$ and $Y$ are continua with $X$ nondegenerate (that is, $X$ is not a point), then $C(X)$ embeds in an ultrapower of $C(Y)$.  Consequently, all nondegenerate models of $T_{cont}$ have the same universal theory.
\end{thm}

We outline the proof here, including extra details not present in the published version communicated to us directly by Hart; we thank him for his permission to include this discussion here.  

By Downward L\"oweinheim-Skolem, we may assume that $C(X)$ and $C(Y)$ are separable, that is, that $X$ and $Y$ are metric continua.  We fix countable bases $\cal B$ and $\cal C$ for their lattices of closed sets and enumerate $\cal C=(C_n)_{n\in \omega}$.  We note that $\cal B^\u$ (the lattice ultrapower) is a base of closed sets for the ultracopower $X_\u$ of $X$.  Thus, by \cite{dowhart}, in order to construct a surjection $X_\u\to Y$ (and thus an embedding $C(Y)\hookrightarrow C(X)^\u$), it suffices to find a map $\phi:\cal C\to \cal B^\u$ satisfying:
\begin{enumerate}
    \item For all $F\in \cal C$, $\phi(F)=\emptyset$ if and only if $F=\emptyset$;
    \item For all $F,G\in \cal C$, if $F\cup G=Y$, then $\phi(F)\cup \phi(G)=(X,X,X,\ldots)_\u$; and
    \item For all $F_1,\ldots,F_n\in \cal C$, if $\bigcap_{i=1}^n F_i=\emptyset$, then $\bigcap_{i=1}^n \phi(F_i)=\emptyset$.
\end{enumerate}

Towards this end, we fix a surjection $f:X\to [0,1]$ (this is where the nondegeneracy of $X$ is used) and identify $Y$ with a closed subspace of the Hilbert cube $Q:= [0,1]^{\bb N}$.  We set $\kappa:\cal P(Y)\to \cal P(Q)$ to be the function $$\kappa(F):=\{q\in Q \ : \ d(q,F)\leq d(q,Y\setminus F)\}.$$ In \cite{kuratowski}, it was shown that this map $\kappa$ has the following properties for all closed $F,G\subseteq Y$:
\begin{itemize}
    \item[(a)] $\kappa(F)\cap Y=F$;
    \item[(b)] $\kappa(F\cup G)=\kappa(F)\cup \kappa(G)$; and
    \item[(c)] $\kappa(Y)=Q$ and $\kappa(\emptyset)=\emptyset$.
\end{itemize}

Moreover, for any $F\subseteq Y$, $\kappa(F)$ is a closed subset of $Q$.  

Set $E:=\{e\subseteq \omega \ : \ \bigcap_{i\in e}C_i=\emptyset\}$ and $E_n:=E\cap \cal P(n)$.  By (a), we have that $Y\cap \bigcap_{i\in e}\kappa(C_i)=\emptyset$ for all $e\in E$.  For each $n\in \omega$, take $\epsilon_n>0$ such that $\epsilon_n<\min\{d(Y,\bigcap_{i\in e}\kappa(C_i)) \ : \ e\in E_n\}$.  Without loss of generality, we may assume that $\lim_{n\to \infty}\epsilon_n=0$.  Since $Y$ is compact, there is a finite open cover $\mathbf{U}_n$ of $U$ by basic open sets contained in the $\epsilon_n$-fattening $Y_{\epsilon_n}:=\{z\in Q \ : \ d(z,Y)\leq \epsilon_n\}$ of $Y$.  Again, since $Y$ is compact, we may fix a finite $\epsilon_n$-dense set $T_n\subseteq Y$ with the property that $T_n\cap C_i\not=\emptyset$ for all $i\leq n$ for which $C_i\not=\emptyset$.  Since $Y$ is connected, any two points in $\bigcup \mathbf{U}_n$ are connected by a piecewise linear path contained in $\bigcup \mathbf{U}_n$.  Consequently, we may define a continuous map $g_n:[0,1]\to \bigcup\mathbf{U}_n\subseteq Q$ containing $T_n$ in its range.  In particular, this map satisfies:
\begin{itemize}
    \item[(i)] $d(g_n(t),Y)<\epsilon_n$ for all $t\in [0,1]$.
    \item[(ii)] For each $y\in Y$, there is $t\in [0,1]$ such that $d(g_n(t),y)<\epsilon_n$.
    \item[(iii)] $g_n^{-1}(C_i)\not=\emptyset$ for all $i\leq n$ for which $C_i\not=\emptyset$.
\end{itemize}

For $i<n<\omega$, set $D_i^n:=f^{-1}(g_n^{-1}(\kappa(C_i)))$, a closed subset of $X$.  By items (b) and (c), whenever $C_i\cup C_j=Y$, we have that $\kappa(C_i)\cup \kappa(C_j)=Q$, whence $D_i^n\cup D_j^n=X$.  Moreover, by (i), whenever $e\in E_n$, we have $g_n^{-1}(\bigcap_{i\in e}\kappa(C_i))=\emptyset$, whence $\bigcap_{i\in e}D^n_i=\emptyset$.  Since $\cal B$ is a lattice base for $X$, there are elements $B_i^n\in \cal B$ containing $D^n_i$ such that $\bigcap_{i\in e}B_i^n=\emptyset$ for all $e\in E_n$.

We are finally ready to define the map $\phi:\cal C\to \cal B^\u$ by setting $\phi(C_i):=(B_i^n)_\u$.  We verify that $\phi$ has the desired properties.  First, if $C_i=\emptyset$, then $\kappa(C_i)=\emptyset$ by (c), whence $D_i^n=\emptyset$ and thus $B_i^n=\emptyset$ (set $e=\{i\}$) for all $n>i$, whence $\phi(C_i)=\emptyset$.  On the other hand, if $C_i\not=\emptyset$, then by (iii), for all $n>i$, we have $g_n^{-1}(C_i)\not=\emptyset$, whence $D_i^n\not=\emptyset$ (since $f$ is surjective) and thus $B_i^n\not=\emptyset$, as desired.

Since $D_i^n\subseteq B_i^n$ for all $i<n<\omega$, it follows from our observations above that whenever $C_i\cup C_j=Y$, we have $B_i^n\cup B_j^n=X$ for all $i<n<\omega$, whence $\phi(C_i)\cup \phi(C_j)=(B_i^n)_\u\cup (B_j^n)_\u=(B_i^n\cup B_j^n)_\u=(X,X,X,\ldots)_\u$, establishing (b).

Finally, suppose that $e\in E_n$.  Then $\bigcap_{i\in e}B^n_i=\emptyset$, whence $$\bigcap_{i\in e}\phi(C_i)=\bigcap_{i\in e}(B_i^n)_\u=\left(\bigcap_{i\in e}B_i^n\right)_\u=\emptyset,$$ establishing (c).  This finishes the proof of Theorem \ref{sameuniversal}.

In a series of papers (see, for example, \cite{bankston}), Bankston studied the model theory of continua in a fairly semantic way by dualizing most notions from classical model theory, occasionally resorting to syntactic techniques by working with lattices bases for the closed sets (which has the disadvantage of not being canonical).  In particular, Bankston introduced the notion of a \textbf{co-existentially closed (co-e.c.) continuum}, which, when formulated in our context, is simply a continuum $X$ for which $C(X)$ is an e.c. model of $T_{cont}$.  Bankston established many properties of co-e.c. continua, including the fact that they are always \textbf{hereditarily indecomposable} (see \cite[Theorem 4.1]{bankston}).  A continuum is indecomposable if it cannot be written as the proper union of two subcontinua and a continuum is hereditarily indecomposable if all subcontinua are indecomposable.  In \cite{bankston2}, Bankston showed that the class of hereditarily indecomposable continua is \textbf{co-elementary}, which, after applying the Gelfand functor, implies that the class of models of $T_{cont}$ of the form $C(X)$ for $X$ a hereditarily indecomposable continuum is elementary.  Moreover, since the inverse limit of hereditarily indecomposable continua is again hereditarily indecomposable, it follows that the aforementioned class of models of $T_{cont}$ is inductive and thus $\forall\exists$-axiomatizable.

One of Bankston's main questions was whether or not the \textbf{pseudo-arc} $\bb P$ is a co-e.c. closed continuum (see \cite[Remark 4.2(i)]{bankston}).  Recall that $\bb P$ is the unique metric continuum that is hereditarily indecomposable and \textbf{chainable}, which is a condition that ensures that the continuum is ``arclike'' in an appropriate sense.  The pseudo-arc is ``generic'' in the descriptive set-theoretic sense of the word \cite{lewis} and thus it is also natural to ask if it is generic from the model-theoretic perspective as well.

In joint work with Eagle and Vignati \cite{EGV}, we were able to answer Bankston's question affirmatively:

\begin{thm}
$C(\bb P)$ is an e.c. model of $T_{cont}$.
\end{thm}

Our main contribution was to show that chainability is a $\forall\bigvee\exists$ property of models of $T_{cont}$ (or rather their image under the inverse Gelfand functor).  By Theorem \ref{sameuniversal}, $C(\bb P)$ is a locally universal model of $T_{cont}$, whence one can enforce being chainable by Proposition \ref{force}.  Since being e.c. is also enforceable, it follows that being e.c. and chainable is enforceable.  Since e.c. implies hereditariliy indecomposable, by the fact that $\bb P$ is the unique metrizable hereditarily indecomposable, chainable continuum, we see that $C(\bb P)$ is actually the enforceable model of $T_{cont}$.

One of the main open questions in the model theory of abelian \cstar-algebras is the following:

\begin{question}
Does $T_{cont}$ have a model companion?
\end{question}

As pointed out in \cite{EGV}, if this model companion exists, then it will \emph{not} be a model completion for $T_{cont}$ does not have the amalgamation property.  Of course, if the model companion exists, then it must in fact be $Th(C(\bb P))$.

% \section{Operator spaces and systems}

% \subsection{The Gurarij Banach space}

% Before embarking on the discussion of the canonical locally universal objects for the classes of operator spaces and systems, it will behoove us to understand the ``classical'' case first, that is, the setting of Banach spaces.  In this case, a canonical locally universal object exists.  In fact, it is the unique separable model of the model completion of the theory of Banach spaces.

% \begin{defn}
% A \textbf{Gurarij space} is a Banach space $X$ having the following property:  for any $\epsilon>0$, any finite-dimensional Banach spaces $E\subseteq F$, and any isometric embedding $T:E\to X$, there is a linear map $S:F\to X$ extending $X$ that satisfies
% $$(1-\epsilon)\|x\|\leq \|S(x)\|\leq (1+\epsilon)\|x\|$$ for all $x\in F$.
% \end{defn}

% In \cite{}, Lusky proved that there is a unique, up to isometric isomorphism, separable Gurarij Banach space, which we shall denote by $\bb G$ and refer to as \textbf{the Gurarij Banach space}.  The uniqueness of $\bb G$ was reproven by Kubis and Solecki in \cite{} and then again by Ben Yaacov in \cite{}, where it was shown that 
% $\bb G$ is the Fra\"isse limit of the class of finite-dimensional...

% \subsection{The noncommutative Gurarij space and the Gurarij operator system}

\end{document}